\documentclass[12pt,reqno]{amsart}
 \usepackage{graphics}
\usepackage{amsfonts, amsmath,enumerate}
\usepackage{amsrefs}
\usepackage{mathtools}

\numberwithin{equation}{section}
\newtheorem{theorem}{Theorem}[section]
\newtheorem{proposition}[theorem]{Proposition}

\newtheorem*{remarks}{Remarks}

\def\n{\nabla}

\def\<{\langle}
\def\>{\rangle}

\def\bea{\begin{eqnarray*} }
\def\eea{\end{eqnarray*} }

\def\beq{\begin{equation}}

\def\R{\mathord{\mathbb R}}

\def\S{\mathord{\mathbb S}}

\def\2{\frac{1}{2}}
\def\3{{\ss}}

\def\.{\cdot}

\def\<{\langle}
\def\>{\rangle}

\def\be{\begin{equation}}
\def\ee{\end{equation}}
\def\bea{\begin{eqnarray}}
\def\eea{\end{eqnarray}}
\def\bsm{\left(\begin{smallmatrix}}
\def\esm{\end{smallmatrix}\right)}
\def\bsmp{\left[\begin{smallmatrix}}
\def\esmp{\end{smallmatrix}\right]}
\def\bpm{\begin{pmatrix}}
\def\epm{\end{pmatrix}}

\newcommand{\quotes}[1]{``#1''}

\begin{document}

\title[Conformally flat isoparametric submanifolds of Euclidean space]{Classification of conformally flat isoparametric submanifolds of Euclidean space}
\author{Christos-Raent Onti}
\address{}
\date{}

\maketitle
\vspace*{-2ex}
\begin{center}
{\small{\it{Dedicated to the memory of Manfredo P. do Carmo}}}
\end{center}

\begin{abstract}
In this note we provide a direct proof of the complete classification of conformally flat 
isoparametric submanifolds of Euclidean space. 
\end{abstract}

\renewcommand{\thefootnote}{\fnsymbol{footnote}} 
\renewcommand{\thefootnote}{\arabic{footnote}} 

\renewcommand{\thefootnote}{\fnsymbol{footnote}} 
\renewcommand{\thefootnote}{\arabic{footnote}}

\section{Introduction}

A submanifold $f\colon M^n\to\R^m$ is said to be {\it isoparametric} if it has flat normal bundle 
and the shape operator in any parallel normal direction has constant eigenvalues. This class of 
Euclidean submanifolds has been investigated extensively by several authors throughout the years; 
see, for example \cite{HPT85,HPT88,BCO16,C94,DFT05,DO04,E01,HL99,KA88,O93,OT75,OT76,PT88,S86,FKM81,T85,T86,T87,T90,TH91,Th00}.
The aim of this short note is to provide a direct proof on the complete classification of conformally flat 
isoparametric submanifolds of Euclidean space by using recent results presented in \cite{DOV18b}. 
Recall that a Riemannian manifold $M^n$ is said to be \emph{conformally flat} if each point 
lies in an open neighborhood conformal to an open subset of Euclidean space 
$\R^n$. 
\vspace{1ex}

The following is the main result.

\begin{theorem}\label{mainconf}
Let $f\colon M^n\to \R^m,\ n\geq 4,$ be an isometric immersion of a conformally flat 
manifold into Euclidean space. If $f$ is isoparametric then $f$ is an open 
subset of one of the following: 
$$
(i)\ \S^n,\ (ii)\ \R^n,\ 
(iii)\ \R^{n-k}\times\S^1\times\dots\times\S^1\subset\R^{n+k},\ 
(iv)\ \ \S^1\times\dots\times\S^1\subset\S^{2n-1}\subset \R^{2n},
$$
$$
\ (v) \ \R\times\S^{n-1}\subset \R^{n+1},\
(vi)\ \S^1\times\S^{n-1}\subset\R^{n+2}.
$$
\end{theorem}

\begin{remarks}
{\rm (I) Theorem \ref{mainconf} can also be extended to non-flat space forms. \\[1mm]
(II) A similar classification also holds if we replace the conformal flatness hypothesis with the
Einstein one, i.e., Riemannian manifolds with constant Ricci curvature. 
In particular, something weaker holds and that is that in this case the submanifold 
does not have to be a priori isoparametric but only to have flat normal bundle and parallel mean curvature vector field
(see \cite{Onti}). The reason we do comment on this case is due to the fact that
these two larger classes, namely the Einstein and the conformally flat ones, are (natural) extensions of the class 
of Riemannian manifolds with constant sectional curvature, which is a subject that has been investigated 
extensively (from the submanifold point of view) by many authors throughout the years; for a survey see \cite{BO01}. 
Notice that the only Riemannian manifolds that are simultaneously Einstein and conformally 
flat are precisely the ones of constant sectional curvature.\\[1mm]
}
\end{remarks}

\vspace{-2ex}

{\noindent {\bf Acknowledgements.} The author would like to thank and express his sincere gratitude 
to Antonio J. Di Scala for his valuable comments and suggestions that have led to the present revised version. 
In particular, he pointed out to the author that: (i) the \quotes{proper} 
assumption (that was in a previous preprint version of Theorem \ref{mainconf}) should be dropped and 
(ii) a different non-direct proof of Theorem \ref{mainconf} can also be derived if one combines deep results 
of Terng, Thorbergsson and Alekseevskii.}

\section{Preliminaries} 

In this section we recall some basic facts. Let $f\colon M^n\rightarrow \R^m$ be an isometric 
immersion of a Riemannian manifold into the Euclidean space.
The second fundamental form $\alpha$ of $f$ is a symmetric 
section of the vector bundle ${\rm Hom}(TM\times TM,N_f M)$, where 
$N_f M$ is the normal bundle of $f$. We say that $f$ is {\it totally umbilical} 
if 
$$
\alpha(X,Y)=\langle X,Y\rangle H,
$$ 
where $H$ is the mean curvature vector field.

If $f$ has flat normal bundle, that is, at any point the curvature tensor 
of the metric induced from the ambient space on the normal bundle of 
the submanifold vanishes, then it is a standard fact (see \cite{RE76}) that at any point 
$x\in M^n$ there exists a set of unique 
pairwise distinct normal vectors $\eta_i(x)\in N_fM(x),\ 1\leq i\leq s(x)$, 
called the {\it principal normals} of $f$ at $x$. 
Moreover, there is an associated orthogonal splitting of the tangent space as
$$
T_xM=E_1(x)\oplus\cdots\oplus E_{s(x)}(x),
$$
where
\begin{equation}\label{E}
E_i(x)=\big\{X\in T_xM:  \alpha(X,Y)=\langle X,Y\rangle\eta_i(x)\ 
\text{for all} \ Y\in T_xM\big\}.
\end{equation}

If $s(x)=k$ is constant on $M^n$, then $f$ is said to be \emph{proper}. 
In this case, the maps 
$
x\in M^n\mapsto\eta_i(x),  \ 1\leq i\leq k,
$ 
are smooth vector fields, called the \emph{principal normal vector fields} of $f$, and the 
distributions 
$
x\in M^n\mapsto E_i(x),\ 1\leq i\leq k,
$ 
are also smooth. If $\n$ denotes the Levi-Civita connection of $M^n$, then 
the Codazzi equation is easily seen to yield 
\begin{equation} \label{cod1} 
\langle\nabla_X Y,Z\rangle(\eta_i-\eta_j) = 
\langle X,Y\rangle\nabla_Z^\perp \eta_i 
\end{equation}
and
\begin{equation}\label{cod3}
\langle\nabla_X V, Z\rangle(\eta_j-\eta_\ell) = 
\langle\nabla_V X, Z\rangle(\eta_j-\eta_i)
\end{equation}
for all $X,Y\in E_i, Z\in E_j,$ and $V\in E_\ell$, where 
$1\leq i\neq j\neq \ell\neq i\leq k$. 
\smallskip

We remark that if $f$ is isoparametric then $f$ is proper.

\vspace{1ex}

The following is contained in \cite{DOV18b}.

\begin{proposition}\label{li} Let $M^n$ be a conformally flat manifold and 
let $f\colon M^n\to\R^m$ be an isometric immersion with flat normal bundle. 
If at some point of $M^n$ we have  $k\geq 3$, then 
the vectors $\eta_j-\eta_i$ and $\eta_j-\eta_\ell$ are linearly 
independent for $1\leq i\neq j\neq \ell\neq i\leq k$.
\end{proposition}

The following is also contained in \cite{DOV18b}.
\vspace{1ex}

\begin{theorem}\label{one} Let $f\colon M^n\to\R^m,n\geq 4$, be an 
isometric immersion with flat normal bundle and  proper of a conformally 
flat manifold. Then $f$ carries at most one 
principal normal vector field of multiplicity larger than one.  
\end{theorem}

The following is well-known; cf. \cite{dt}.

\begin{proposition}\label{conflat}
A Riemannian product is conformally flat if and only if one of the following 
possibilities holds:
\begin{enumerate}[(i)]
  \item One of the factors is one-dimensional and the other one has constant sectional curvature.
  \item Both factors have dimension greater than one and are either both flat or have opposite 
  constant sectional curvatures.
\end{enumerate}
\end{proposition}

A map $f\colon M^n\rightarrow \R^m$ from a product manifold 
$M^n=\Pi_{i=1}^k M_i$ is called the {\it extrinsic product of immersions} 
$
f_i\colon M_i\rightarrow \R^{m_i},\ 1\leq i\leq k,
$ 
if there exist 
an orthogonal decomposition $\R^m=\Pi_{i=0}^k \R^{m_i}$, with 
$\R^{m_0}$ possibly trivial, such that $f$ is given by 
$$f(x)=(v,f_1(x_1),\dots,f_k(x_k))$$ for all $x=(x_1,\dots,x_k)\in M^n$ 
and $v\in\R^{m_0}$.

Let $f\colon M^n \to\R^m$ be an isometric immersion of a
Riemannian manifold. If $M^n=\Pi_{i=1}^k M_i$ is a product manifold 
then the second fundamental form $\alpha$ is  said to be {\it adapted} to the product structure of $M^n$ if 
$$\alpha(X_i,X_j)=0\ \text{for all} \ X_i\in TM_i, \ X_j\in TM_j\ \text{with}\ 1\leq i\neq j\leq k,$$ 
where the tangent bundles $TM_i$ are identified with the corresponding tangent distributions
to $M^n$.
The next result, due to Moore \cite{MO71}, shows that extrinsic products of isometric immersions are 
characterized by this property among isometric immersions of Riemannian products.

\begin{theorem}\label{Moore}
Let $f\colon M^n\rightarrow \R^m$ be an isometric immersion of a Riemannian product manifold $M^n=\Pi_{i=1}^k M_i$ 
with adapted second fundamental form. Then $f$ is an extrinsic product of isometric immersions.
\end{theorem}

\section{Proof of Theorem \ref{mainconf}}

Assume that $k\geq 3$, since 
otherwise the result is immediate. We claim that each distribution 
$E_i,\ 1\leq i\leq k,$ is parallel, that is 
$$\nabla_X Y\in E_i\ \text{for all}\ X\in TM,\ Y\in E_i \ \text{and}\ 1\leq i\leq k.$$

Since $f$ is isoparametric we have that the principal normal 
vector fields are parallel in the normal connection. Therefore, it follows from
the Codazzi equation \eqref{cod1} that each distribution $E_i,\ 1\leq i\leq k,$ is totally 
geodesic. Thus, we only need to show that 
$$
\nabla_X Y\in E_i\ \text{for all}\ X\in E_j\ \text{and}\ Y\in E_i
$$ 
with $j\neq i$. Indeed, we consider $Y\in E_i,\ X\in E_j,\ Z\in E_\ell\subset E_i^\perp$ 
and distinguish the following two cases. 

If $\ell=j$, then we get
\begin{equation}\label{par1}
\langle \nabla_X Y,Z\rangle=-\langle Y,\nabla_X Z\rangle=0,
\end{equation}
where we have used the fact that $E_\ell$ is totally geodesic.

If $\ell\neq j$, then from \eqref{cod3} we obtain
$$\langle \nabla_X Y,Z\rangle (\eta_\ell-\eta_i)=\langle \nabla_Y X,Z\rangle (\eta_\ell-\eta_j).$$
Using Proposition \ref{li}, we get  
\begin{equation}\label{par2}
\langle \nabla_X Y,Z\rangle=0 \ \text{for all}\ X\in E_j, Y\in E_i, Z\in E_\ell, 
\end{equation}
with $\ell\neq i\neq j$.
Therefore, from \eqref{par1} and \eqref{par2}, we obtain that 
$\nabla_X Y\in E_i$ for all $X\in E_j$ and $Y\in E_i$ with $j\neq i$. 
This completes the proof of the claim.

Now, de Rham's theorem implies that around 
every point $x\in M^n$ there is a neighborhood $U$ that is the Riemannian product of the integral manifolds 
$M_1,\dots,M_k$ of the distributions $E_1,\dots,E_k$ respectively, through a point $y\in U$. 
Therefore, since the second fundamental form of $f$ is adapted, 
Theorem \ref{Moore} implies that $\left. f\right\vert_U$ is an extrinsic product of 
isometric immersions $f_i\colon M_i\to \R^{m_i},\ 1\leq i\leq k$, which due to \eqref{E} and our hypothesis 
have to be totally umbilical and with mean curvatures of constant length. 
Now, the classification follows easily by using Theorem \ref{one} and Proposition \ref{conflat}. 
This completes the proof. \qed

\end{document}